# CORRECTION
# PERFECT SIMULATION FOR A CLASS OF POSITIVE RECURRENT MARKOV CHAINS
*Ann. Appl. Probab.* **17** (2007) 781–808

By Stephen B. Connor and Wilfrid S. Kendall

*University of Warwick*

In [1] we introduced a class of positive recurrent Markov chains, named tame chains. A perfect simulation algorithm, based on the method of dominated CFTP, was then shown to exist in principle for such chains. The construction of a suitable dominating process was flawed, in that it relied on an incorrectly stated lemma ([1], Lemma 6). This claimed that a geometrically ergodic chain, subsampled at a stopping time $\sigma$, satisfies a geometric Foster–Lyapunov drift condition with coefficients not depending on $\sigma$. This is true if $\sigma$ is a stopping time independent of the chain, but *not* if this independence does not hold. Reference [1], Lemma 6 is therefore false as stated.

We now indicate a corrected construction of a dominating process. As described in [1], Section 3.1, the process $D$ is defined by starting with a process $Y$ and pausing it using a function $S$. In the following modified construction this is simplified by taking $S = F$, where $F$ is the function taming $X$. We restate [1], Theorem 16, and give a shorter proof, which avoids the faulty Lemma 6 but pays a price in terms of consequences for the perfect simulation algorithm of Section 3.3. The discussion of tameness (Section 4) is unaffected.

THEOREM 16. *Suppose $X$ satisfies the weak drift condition $PV \leq V + b\mathbf{1}_C$, and that $X$ is tamed with respect to $V$ by the function*

$$F(z) = \begin{cases} \lceil \lambda z^\delta \rceil, & z > d', \\ 1, & z \leq d', \end{cases}$$

*with the resulting subsampled chain $X'$ satisfying a drift condition $PV \leq \beta V + b'\mathbf{1}_{[V \leq d']}$, with $\log \beta < \delta^{-1} \log(1-\delta)$. Then there exists a stationary ergodic process $D$ which dominates $V(X)$ at the times $\{\sigma_n\}$ when $D$ moves.*

PROOF. Suppose that $D_{\sigma_n} = z$, and that $V(X_{\sigma_n}) = V(x) \leq z$. We wish to show that $D_{\sigma_{n+1}}$ can dominate $V(X_{\sigma_{n+1}})$, where $\sigma_{n+1} = \sigma_n + F(z)$ is the time at which $D$ next moves. Domination at successive times $\sigma_j$ at which $D$ moves then follows inductively. For simplicity in the calculations below we set $\sigma_n = 0$.

First choose $\beta^* > \beta$ such that

(1) $$\log \beta < \log \beta^* < \delta^{-1} \log(1-\delta).$$

Our aim is to control $\mathbb{E}_x[V(X_{F(z)})]$, recalling that $F(z)$ is deterministic and that $F(V(x)) \leq F(z)$:

$$\mathbb{E}_x[V(X_{F(z)})] = \mathbb{E}_x[V(X_{F(V(x))})] + \mathbb{E}_x[V(X_{F(z)}) - V(X_{F(V(x))})]$$
$$= \mathbb{E}_x[V(X'_1)] + \mathbb{E}_x[V(X_{F(z)}) - V(X_{F(V(x))})]$$
$$\leq \beta V(x) + b' \mathbf{1}_{[V(x) \leq d']} + b[F(z) - F(V(x))]$$
$$\leq \beta z + b' + b(\lambda + 1) z^\delta$$

(2) $$\leq \beta^* z \quad \text{for } z \geq h^*,$$

where $h^* < \infty$ is a constant chosen sufficiently large for inequality (2) to hold. The first inequality in this sequence holds due to the drift conditions satisfied by $X'$ and $X$. The second follows from the definition of $F$ and the assumption that $V(x) \leq z$.

Now define the process $Y = h^* \exp(U)$, where $U$ is the system workload of a $D/M/1$ queue with arrivals every $\log(1/\beta^*)$ time units and service times being independent and of unit Exponential distribution. As in the original proof of Theorem 16, $Y$ may be paused using $F$ to obtain the process $D$ which is positive recurrent and has a proper equilibrium distribution by virtue of inequality (1).

Finally, observe that $D$ takes values in $[h^*, \infty)$. As in the proof of Theorem 5 of [2], it follows from inequality (2) that $V(X_{F(z)})$ can be dominated by $D_{F(z)}$, as required. $\square$

The majority of Section 3.3 remains valid when the dominating process is constructed as above. The only issue is that by taking $S = F$ in this new method we are no longer assured that $S(h^*) \geq m$, where the set $C^* = \{x : V(x) \leq h^*\}$ is $m$-small. Unfortunately, there no longer seems to be a simple way to ensure this since our attempts to increase $S$ in the above always result in an increased value of $h^*$.

If it so happens that $F(h^*) \geq m$ for a given chain, then the original perfect simulation algorithm remains unchanged. If this is not the case, then the algorithm must be altered. It now becomes necessary, when $D_0 = h^*$, for $D$ to dominate $V(X)$ not at time $\sigma_1 = F(h^*)$ but at time

$$\sigma^* = \inf_{j \geq 2} \{\sigma_j : \sigma_j \geq m\}.$$

This is an example of the composite nondeterministic sampling schemes we had originally hoped to avoid (cf. the comment before [1], Theorem 15]). Furthermore, we need to be able to couple target chains and dominating process at $\sigma^*$ in such a way that the target chains may regenerate at this time (using the fact that $C^*$ is $\sigma^*$-small). This unfortunately reduces the impact of the result, which is an issue that we are currently trying to resolve.

DEPARTMENT OF STATISTICS
UNIVERSITY OF WARWICK
COVENTRY CV4 7AL
UNITED KINGDOM
E-MAIL: s.b.connor@warwick.ac.uk
       w.s.kendall@warwick.ac.uk